\documentclass{amsart}

\usepackage{xypic}
\input xy
\xyoption{all}
\usepackage{epsfig}
\usepackage{amsthm}
\usepackage{amssymb}
\usepackage{amsmath}
\usepackage{amscd}
\usepackage{color}
\usepackage[T1]{fontenc}

\newcommand{\bg}{\begin{equation}}
\newcommand{\ed}{\end{equation}}
\newcommand{\bga}{\begin{eqnarray}}
\newcommand{\eda}{\end{eqnarray}}
\newcommand{\pf}{\textbf{Proof:\ }}

\def\cbdu{\par{\raggedleft$\Box$\par}}

\newtheorem {Theorem}  {Theorem}

\numberwithin{Theorem}{section}

\newtheorem {Lemma}[Theorem]  {Lemma}

\theoremstyle{definition}
\newtheorem{Definition}[Theorem]{Definition}
\theoremstyle{remark}

\newtheorem {Corollary}[Theorem]{\bf Corollary}

\expandafter\chardef\csname pre amssym.def
at\endcsname=\the\catcode`\@ \catcode`\@=11
\def\undefine#1{\let#1\undefined}
\def\newsymbol#1#2#3#4#5{\let\next@\relax
 \ifnum#2=\@ne\let\next@\msafam@\else
 \ifnum#2=\tw@\let\next@\msbfam@\fi\fi
 \mathchardef#1="#3\next@#4#5}
\def\mathhexbox@#1#2#3{\relax
 \ifmmode\mathpalette{}{\m@th\mathchar"#1#2#3}%
 \else\leavevmode\hbox{$\m@th\mathchar"#1#2#3$}\fi}
\def\hexnumber@#1{\ifcase#1 0\or 1\or 2\or 3\or 4\or 5\or 6\or 7\or 8\or
 9\or A\or B\or C\or D\or E\or F\fi}

\font\teneufm=eufm10 \font\seveneufm=eufm7 \font\fiveeufm=eufm5
\newfam\eufmfam
\textfont\eufmfam=\teneufm \scriptfont\eufmfam=\seveneufm
\scriptscriptfont\eufmfam=\fiveeufm

\catcode`\@=\csname pre amssym.def at\endcsname

\newcounter{remark}
\newcommand{\vf}{\varphi}

\newcommand{\R}{\mathbf{R}}

\newcommand{\Ll}{{\mathcal L}}

\def  \I   {{\mathbb I}}

\def  \R   {{\mathbb R}}

\def  \Q   {{\mathbb Q}}
\def  \S   {{\mathbb S}}

\def  \H   {{\mathbb H}}
\def  \A   {{\mathbb A}}
\def  \B   {{\mathbb B}}

\def  \12  {{\frac{1}{2}}}

\def\build#1_#2^#3{\mathrel{\mathop{\kern 0pt#1}\limits_{#2}^{#3}}}

\begin{document}
%\currannalsline{0}{2006}

\title[Regularity for LCD with Q-tensor]{Regularity problem for the nematic LCD system with Q-tensor in $\mathbb R^3$}

%\author{hello}

\author [Mimi Dai]{Mimi Dai}
\address{Department of Mathematics, Stat. and Comp. Sci.,  University of Illinois Chicago, Chicago, IL 60607,USA}
\email{mdai@uic.edu} 

%\thanks{The work of the authors was partially supported by NSF Grant DMS--1108864.}

%%%use \Proof instead of \begin{proof}

%%%% use \Endproof instead of \end{proof}

%%%% use \references {999} instead of \begin{thebibliography}{99}

%%%%used \Endrefs instead of \end{thebibliography}

\begin{abstract}
We study the regularity problem of a nematic liquid crystal model with local configuration represented by Q-tensor in three dimensions. It was an open question whether the classical Prodi-Serrin condition implies regularity for this model. Applying a wavenumber splitting method, we show that a solution does not blow-up under certain extended Beale-Kato-Majda condition solely imposed on velocity. This regularity criterion automatically implies that the classical Prodi-Serrin or Beale-Kato-Majda condition prevents blow-up of solutions.

\bigskip


\hspace{0.02cm}CLASSIFICATION CODE: 35B44, 35Q35, 76A15, 76D03.
\end{abstract}

\maketitle

\section{Introduction}

Considered here is a hydrodynamic model of nematic liquid crystals proposed by Beris and Edwards \cite{BE}, where the local configuration of the crystal is represented by the $Q-$tensor $\Q=\Q(x,t)$.  The evolution of the crystal flow is governed by 
\begin{equation}\label{LCD}
\begin{split}
u_t+(u\cdot\nabla) u+\nabla p=\nu\Delta u+\nabla\cdot \Sigma(\Q),\\
\Q_t+(u\cdot\nabla )\Q-\S(\nabla u,\Q)=\mu\Delta\Q-\Ll[\partial F(\Q)]\\
\nabla \cdot u=0,
\end{split}
\end{equation}
with $(x,t)\in \R^n\times(0,\infty)$.
In the equations,  the $Q-$tensor $\Q$ describes the ordering of the molecule, $u$ is the fluid velocity, $p$ is the fluid pressure and $F$ denotes a potential function which will be described later. The parameter $\nu$ denotes the kinematic viscosity coefficient of the fluid, and $\mu$ stands for the elasticity of the molecular orientation field. 
The tensor $\Q\in \R^{3\times 3}_{{\rm sym},0}$ is symmetric and traceless. Other notations are introduced in the following. The operator $\Ll$ is defined as
\[\Ll[\A]=\A-\frac13{\rm tr}[\A]\I\]
which represents the projection onto the space of traceless matrices. The tensors $\Sigma$ and $\S$ are given by
\begin{equation}\label{sig}
\Sigma(\Q)=2\xi \H:\Q \left(\Q+\frac13\I\right)-\xi \left[\H \left(\Q+\frac13\I\right)-\left(\Q+\frac13\I\right)\H\right]
-(\Q\H-\H\Q)-\nabla\Q\otimes\nabla \Q,
\end{equation}
\begin{equation}\label{S}
\S(\nabla u,\Q)=(\xi D(u)+\Omega(u))\left(\Q+\frac13\I\right)+\left(\Q+\frac13\I\right)(\xi D(u)-\Omega(u))
-2\xi \left(\Q+\frac13\I\right)\Q:\nabla u,
\end{equation}
with $(\nabla\Q\otimes\nabla \Q)_{ij}=\partial_i\Q_{\alpha\beta}\partial_j\Q_{\alpha\beta}$. Note that we use the summation convention for repeated indices here and through the rest of the paper. 
In the above equations, $D$ and $\Omega$ denote the symmetric and skew-symmetric parts of the rate of strain tensor, respectively
\begin{equation}\notag
D(u)=\frac12(\nabla u+\nabla^tu), \ \quad \Omega(u)=\frac12(\nabla u-\nabla^tu).
\end{equation}
While $\H$ is obtained through the variational derivative of the free energy under the constraint that $\Q$ is symmetric and traceless, as
\[\H=\Delta \Q-\Ll[\partial F(\Q)],\]
where $F(\Q)$ denotes the bulk potential function. In this work, we take the Landau-de Gennes form
\begin{equation}\label{F}
F(\Q)=\frac a2|\Q|^2+\frac b3 \rm{tr}[\Q^3]+\frac c4|\Q|^4
\end{equation}
with constants $a, b, c\in \R$ which are determined by the material and temperature. More general choices of $F(\Q)$ are considered in \cite{BM, DFRSS, Wil}. 
The coefficient $\xi\in\R$
measures the ratio between the rotation and aligning effects that a shear flow exerts over the directors.  We assume $\xi=0$ for simplicity. Thus, (\ref{sig}) and (\ref{S}) become
\begin{equation}\label{tensors}
 \Sigma(\Q)=\Delta \Q\Q-\Q\Delta\Q-\nabla \Q\otimes\nabla\Q, \ \quad   \S(\nabla u,\Q)=\Omega(u)\Q-\Q\Omega(u)
\end{equation}
after using 
\[\Q\Ll[\partial F(\Q)]-\Ll[\partial F(\Q)]\Q=\Q\partial F(\Q)-\partial F(\Q)\Q=0.\]
We expect the same result shall hold for the general case $\xi\neq0$, which will be considered in a future work.
For more discussions regarding the physics of the model, we direct the readers to \cite{BE} and the references therein.

Regarding the mathematical aspects of the model, we briefly mention a few fundamental and relevant results in the literature without the intention to be complete. The existence of weak solutions to (\ref{LCD}) was established by Paicu and Zarnescu \cite{PZenergy, PZex} in both two and three dimensions (3D), for $\xi=0$ and small $|\xi|>0$ respectively. Moreover, in 2D, the authors obtained global regular solutions. %{\color{blue}Existence for general $\xi$ in 2D in \cite{CRWX}} 
On a bounded domain in 3D, Abels et al.  \cite{ADL, ADL13} proved the existence and uniqueness of local strong solution subject to various boundary conditions. Regarding the long time behavior,  Dai et al.  \cite{DFRSS} established the optimal decay rate for {\it weak} solutions in 3D, while optimal decay rates are usually obtained for regular solutions for other liquid crystal models (see \cite{DQS, DS, Wu}). Another interesting work is by Cavaterra et al. \cite{CRWX} who showed the existence of global strong solutions and decay rate of strong solutions in 2D with general $\xi\in\R$. 

In this paper, we study the global regularity problem for the Q-tensor model (\ref{LCD}) in 3D. 
The existence of global regular solutions for the 3D Navier-Stokes equation (NSE) is an outstanding open problem, and thus an open problem for (\ref{LCD}) as well. 
For the 3D NSE, it is well known that the Prodi-Serrin condition 
\begin{equation}\label{PS}
u\in L^l(0,\infty; L^r), \qquad \frac{2}{l}+\frac{3}{r}=1, \quad r>3
\end{equation}
or the Beale-Kato-Majda (BKM) condition (see \cite{BKM})
\begin{equation}\label{BKM-criterion}
\int_0^T\|\nabla \times u\|_\infty\, dt<\infty,
\end{equation}
guarantees global regularity.  For the invisid NSE (Euler equation), Planchon \cite{P} gave an extended BKM  condition for regularity, namely
\begin{equation}\label{BKM-extended}
\lim_{\epsilon\to 0}\sup_{q}\int_{T-\epsilon}^T\|\Delta_q(\nabla\times u)\|_{\infty}\, dt<c,
\end{equation}
for a small enough constant $c$, where $\Delta_q$ denotes the Littlewood-Paley projection (see Section \ref{sec:pre}). Another improvement of the BKM criterion is given by Cheskidov and Shyvdkoy \cite{CS}, as
\begin{equation} \label{CS-criterion}
\int_0^T\|\nabla \times u_{\leq Q}\|_{B^0_{\infty,\infty}}\, dt<\infty, \qquad \mbox {for some wavenumber} \qquad \Lambda(t)=2^{Q(t)},
\end{equation}
where $u_{\leq Q}$ denotes the low modes part of the velocity (see Section \ref{sec:pre}).
Condition (\ref{CS-criterion}) is also weaker than all the Prodi-Serrin criteria (\ref{PS}). A wavenumber splitting method was introduced in \cite{CS} to achieve the goal. Recently, Cheskidov and Dai \cite{CD} refined the method and established the following regularity criterion 
\bg\label{CD-criterion}
\limsup_{q \to \infty} \int_{T/2}^T 1_{q\leq Q(\tau)}\lambda_q\|u_q\|_\infty \, d\tau  \leq c,
\ed
with a small constant $c$, and certain wavenumber $\Lambda(t)=2^{Q(t)}$. This condition is weaker than all the conditions above, (\ref{PS}), (\ref{BKM-criterion}), (\ref{BKM-extended}) and (\ref{CS-criterion}). 
%Indeed, \[u_{\leq Q} \in L^r(0,T; B^{-1+\frac{2}{r}}_{\infty,\infty}) \qquad \text{for some} \qquad 1\leq r <\infty\] implies \eqref{eq:BKM-Lambda}. Note that $r$ cannot be taken as $\infty$ here. 
%A regularity criterion in the limit case $r=\infty$ has been established by Escauriaza, Seregin, and Sverak \cite{ESS} in $L^\infty(0;T,L^3)$, and extended by Gallagher, Koch, and Planchon to $L^\infty(0,T;B^{-1+ \frac{3}{s}}_{s,q})$ with $3 < s,q< \infty$ in \cite{GKP}.
Back to the Q-tensor model, it seems that not much previous work on this topic can be found in the literature, though Guill\' en-Gonz\'alez and Rodr\'iguez-Bellido pointed out that the classical Prodi-Serrin regularity criteria are not valid for the model due to the presence of the stretching terms (c.f. Remark 4, \cite{GGRB15}). In another work of the same authors \cite{GGRB}, they studied the so-called {\it weak regularity} for $(\partial_tu,\partial_t\Q)$ which is different from the standard regularity problem as in the Navier-Stokes equation framework. They actually need to impose conditions on both of the velocity and the Q-tensor to imply the global in time of  {\it weak regularity}. Regarding the standard (strong) regularity, they showed that with an additional assumption on the gradient of velocity as
\begin{equation}\notag
\nabla u\in L^{p}(0,T; L^q), \ \quad \frac2p+\frac3q=2, \ \quad 2\leq q\leq 3,
\end{equation} 
the solution on $(0,T)$ does not blow-up at time $T$.
In this paper, we shall apply the wavenumber splitting method to the Q-tensor model and obtain that a condition analogous to (\ref{CS-criterion}) or (\ref{CD-criterion}) solely imposed on velocity yields global regularity of (\ref{LCD}) in 3D. Thus, automatically, the classical Prodi-Serrin or BKM condition solely on velocity is sufficient to imply global regularity for the model.

Based on the scaling of the flow equation in (\ref{LCD}), we define the dissipation wavenumber $\Lambda(t)$ for the velocity as
\begin{equation}\label{wave}
\Lambda(t)=\min \left\{\lambda_q:\lambda_p^{-1+\frac 3r}\|u_p(t)\|_r<c_r\min\{\nu,\mu\}, \forall p>q, q\in \mathbb{N} \right\},
\end{equation}
with $\lambda_q =2^q$. For more notations, $u_p = \Delta_p u$ denotes the Littlewood-Paley projection of $u$,  and $c_r$ is a dimensionless  constant that depends only on $r\in[2,6)$. Take $Q(t)\in\mathbb N$ such that $\lambda_{Q(t)}=\Lambda(t)$.
Our main results are stated below.
\begin{Theorem}\label{thm}
Let $(u, \Q)$ be a weak solution to (\ref{LCD}) on $[0,T]$. Assume that
\bg\label{criterion}
\int_0^T \|\nabla u_{\leq Q}(t)\|_{B^0_{\infty,\infty}} \, dt  \leq \infty.
%\limsup_{q \to \infty} \int_{T/2}^T 1_{q\leq Q(\tau)}\lambda_q\|u_q\|_\infty \, d\tau  \leq c_r,
\ed
%with some $2\leq r \leq \infty$.
Then $(u(t), \Q(t))$ is regular on $(0,T]$.
\end{Theorem}

\begin{Theorem}\label{thm-strong}[A stronger statement.]
Let $(u, \Q)$ be a weak solution to (\ref{LCD}) on $[0,T]$. Assume that $(u(t), \Q(t))$ is regular on $(0,T)$, and
\bg\label{criterion2}
\limsup_{q \to \infty} \int_{T/2}^T 1_{q\leq Q(\tau)}\lambda_q\|u_q\|_\infty \, d\tau  \leq c,
\ed
for a small constant $c$. Then $(u(t), \Q(t))$ is regular on $(0,T]$.
\end{Theorem}

\begin{Corollary}\label{cor}[Prodi-Serrin and BKM type criteria.] 
Let $(u, \Q)$ be a weak solution to (\ref{LCD}) as in Theorem \ref{thm}. Assume that one of the following holds
\bg\label{criterion3}
u\in L^s(0,T;L^r(\R^3)), \ \quad \mbox { with } \frac2s+\frac3r=1, \ \quad 3<r<6;
\ed
\bg\label{criterion4}
\int_0^T \|\nabla \times u(t)\|_{\infty} \, dt  \leq \infty.
\ed
Then $(u(t), \Q(t))$ is regular on $(0,T]$.
\end{Corollary}

The results have certainly a novelty value. First of all, each condition is solely imposed on the low frequency part of the
velocity despite the fact that (\ref{LCD}) is a coupled system of the NSE and the evolution of the Q-tensor. Second, even though strong nonlinearity appears in the system, for instance the term $\nabla\cdot (\H \Q-\Q\H)$, we are able to deal with it by revealing cancelations. More over, as stated in Corollary \ref{cor}, the classical Prodi-Serrin and BKM type criteria are valid for the Q-tensor system, which was previously unknown.  Although, we would like to point out that there is a restriction on the index $r<6$ for the Prodi-Serrin type criteria obtained in Corollary \ref{cor}. Extending it to the full range $r\in(3,\infty]$ will be addressed in a future work. 

We mention that the wavenumber splitting method has also been successfully applied to other disipative equations, for instance, the supercritical surface quasi-geostrophic equation in \cite{Dsqg}, the magneto-hydrodynamics system in \cite{CD}, and the Hall magneto-hydrodynamics system in \cite{Dhmhd}.

The rest of the paper is planned as follows. In Section \ref{sec:pre} we introduce some notations, recall the Littlewood-Paley theory, describe the energy laws of solutions to (\ref{LCD}), and establish some commutator estimates. Section \ref{sec:reg} is devoted to proving Theorem \ref{thm}. The proof of Theorem \ref{thm-strong} will be omitted since an identical analysis based on the proof of Theorem \ref{thm} can be found in previous work \cite{CD}. On the other hand, (\ref{criterion4}) obviously implies (\ref{criterion}); and it is shown in \cite{CD} that (\ref{criterion3}) implies (\ref{criterion}) as well. Thus, the proof of Corollary \ref{cor} will be not present either.

\bigskip

\section{Preliminaries}
\label{sec:pre}

\subsection{Notation}
\label{sec:notation}
We denote by $A\lesssim B$ an estimate of the form $A\leq C B$ with
some absolute constant $C$, and by $A\sim B$ an estimate of the form $C_1
B\leq A\leq C_2 B$ with some absolute constants $C_1$, $C_2$. 
We write $\|\cdot\|_p=\|\cdot\|_{L^p}$. The symbol $(\cdot, \cdot)$ stands for the $L^2$-inner product.

Let $S_0^{d}$ be the space of $Q$-tensors in dimension $d$, namely,
\[S_0^d=\{\Q\in \mathbb M^{d\times d}: \Q_{ij}=\Q_{ji}, tr(\Q)=0, i,j=1, ..., d\}.\]
We adapt the Frobenius norm of a matrix $|\Q|=\sqrt{tr \Q^2}=\sqrt{\Q_{\alpha\beta}\Q_{\alpha\beta}}$ and define Sobolev spaces of $Q$-tensors in terms of this norm.

\subsection{Littlewood-Paley decomposition}
\label{sec:LPD}
Our method  relies on the Littlewood-Paley decomposition, which we briefly recall here. 
For a more detailed description on this theory, we direct the readers to the books \cite{BCD} and \cite{Gr}.

We denote the Fourier transform and inverse Fourier transform by $\mathcal F$ and $\mathcal F^{-1}$, respectively.  Let $\chi\in C_0^\infty(\R^n)$ be a nonnegative radial function such that 
\begin{equation}\notag
\chi(\xi)=
\begin{cases}
1, \ \ \mbox { for } |\xi|\leq\frac{3}{4}\\
0, \ \ \mbox { for } |\xi|\geq 1.
\end{cases}
\end{equation}
More bump functions are chosen as
\begin{equation}\notag
\varphi(\xi)=\chi(\xi/2)-\chi(\xi), \ \quad
\varphi_q(\xi)=
\begin{cases}
\varphi(\lambda_q^{-1}\xi)  \ \ \ \mbox { for } q\geq 0,\\
\chi(\xi) \ \ \ \mbox { for } q=-1,
\end{cases}
\end{equation}
with $\lambda_q=2^q$. Note that the sequence of the smooth functions $\{\varphi_q\}_{q\geq-1}$ forms a dyadic partition of unity.
For a tempered distribution vector field $u$ we define the Littlewood-Paley decomposition
\begin{equation}\notag%\label{eq:LPP}
\begin{cases}
&h=\mathcal F^{-1}\varphi, \qquad \tilde h=\mathcal F^{-1}\chi,\\
&u_q:=\Delta_qu=\mathcal F^{-1}(\varphi(\lambda_q^{-1}\xi)\mathcal Fu)=\lambda_q^n\int h(\lambda_qy)u(x-y)dy,  \qquad \mbox { for }  q\geq 0,\\
& u_{-1}=\mathcal F^{-1}(\chi(\xi)\mathcal Fu)=\int \tilde h(y)u(x-y)dy.
\end{cases}
\end{equation}
Then
\bg\notag
u=\sum_{q=-1}^\infty u_q
\ed
holds in the distributional sense.  To simplify the notation, we denote
\bg\notag
\tilde u_q=u_{q-1}+u_q+u_{q+1}, \qquad u_{\leq Q}=\sum_{q=-1}^Qu_q,  \qquad u_{(P,Q]}=\sum_{q=P+1}^Qu_q.
\ed

The Besov space $B_{p,\infty}^{s}$ is defined as follows.
\begin{Definition}
Let $s\in \mathbb R$, and $1\leq p\leq \infty$. The Besov space $B_{p,\infty}^{s}$ is the space of tempered distributions $u$ such that the following norm
$$
\|u\|_{B_{p, \infty}^{s}}=\sup_{q\geq -1}\lambda_q^s\|u_q\|_p
$$
is finite.
\end{Definition}
Note that 
\[
  \|u\|_{H^s} \sim \left(\sum_{q=-1}^\infty\lambda_q^{2s}\|u_q\|_2^2\right)^{1/2},
\]
for each $u \in H^s$ and $s\in\R$.

We recall Bernstein's inequality (c.f. \cite{L}).
\begin{Lemma}\label{le:bern} %\cite{L}
Let $n$ be the space dimension and $r\geq s\geq 1$. Then for all tempered distributions $u$, 
\bg\notag%\label{Bern}
\|u_q\|_{r}\lesssim \lambda_q^{n(\frac{1}{s}-\frac{1}{r})}\|u_q\|_{s}.
\ed
\end{Lemma}

\medskip

\subsection{Energy law, weak and strong solutions, maximum principle}
\label{sec:sol}

Denote the energy functional as 
\[E=\frac12|u|^2+\frac12|\nabla \Q|^2+F(\Q).\]
The basic energy law is given by (c.f. \cite{PZenergy})
\[\frac {d}{dt}\int_{\R^3}E\, dx+\int_{\R^3}|\nabla u|^2+\left|\Delta \Q-\Ll[\partial F(\Q)]\right|^2\, dx=0\]
provided that $u$ and $\Q$ vanish for large $|x|$.

We recall the standard definitions of weak and regular solutions to a differential equation system (see, e.g.,  \cite{Tem}).
\begin{Definition}\label{def:weak} 
A weak solution of (\ref{LCD}) on $[0,T]$ is a pair of functions $(u, \Q)$ in the class 
\begin{equation}\notag
\begin{split}
&u\in C_w([0,T]; L^2(\mathbb R^3)) \cap L^2(0,T; H^1(\mathbb R^3)),  \\
&\Q\in C_w([0,T]; H^1(\mathbb R^3)) \cap L^2(0,T; H^2(\mathbb R^3)),
\end{split}
\end{equation}
satisfying 
\begin{equation}\notag%\label{eq:uweak}
\begin{split}
&(u(t), \vf(t))-(u_0, \vf(0))-\int_0^t(u(s), \partial_s\vf(s))+\nu(u(s), \Delta\phi(s))\, ds\\
=&\int_0^t (u(s)\cdot\nabla\vf(s), u(s))+([\Q(s)\H(s)-\H(s)\Q(s)]+\nabla\Q(s)\otimes\nabla\Q(s), \nabla\vf(s))\, ds,
\end{split}
\end{equation}
\begin{equation}\notag%\label{eq:bweak}
\begin{split}
&(\Q(t), \psi(t))-(\Q_0, \psi(0))-\int_0^t (\Q(s), \partial_s\psi(s))+\mu(\Q(s), \Delta\psi(s))\, ds\\
=&\int_0^t (u(s)\cdot\nabla\psi(s), \Q(s))-\left(\Omega(u(s))\Q(s)-\Q(s)\Omega(u(s)), \psi(s)\right)\, ds,
\end{split}
\end{equation}
for all smooth functions $\vf\in C_0^\infty([0,T]\times\mathbb R^3; \mathbb R^3)$ with $\nabla_x\cdot \vf=0$, and $\psi\in C_0^\infty([0,T]\times\mathbb R^3; S_0^3)$.
\end{Definition}

\begin{Definition}\label{le-reg} 
A weak solution $(u, \Q)$ of (\ref{LCD}) is regular on a time interval $\mathcal I$ if $\|u(t)\|_{H^s}$ and $\|\nabla\Q(t)\|_{H^s}$ are continuous on $\mathcal I$ for some $s> \frac{1}{2}$.
\end{Definition}

In \cite{DFRSS}, we proved the following maximum principle for the $Q$-tensor equation.
\begin{Lemma}
\label{max}
Let $(u,\Q)$ be a weak solution to (\ref{LCD}) with initial data $(u_0,\Q_0)$. Then, for all $t>0$, it holds
\begin{equation}\notag
\|\Q(\cdot,t)\|_{L^\infty(\R^3;\R^{3\times3})}\leq  \|\Q_0\|_{L^\infty(\R^3;\R^{3\times3})}.
\end{equation}
\end{Lemma}

\subsection{Commutators}

To deal with the Littlewood-Paley projection of a product term, we often decompose it by the so-called Bony's paraproduct according to different types of interactions. Namely, for instance, we write
\begin{equation}\notag
\begin{split}
\Delta_q(u\cdot\nabla v)=&\sum_{|q-p|\leq 2}\Delta_q(u_{\leq{p-2}}\cdot\nabla v_p)+
\sum_{|q-p|\leq 2}\Delta_q(u_{p}\cdot\nabla v_{\leq{p-2}})\\
&+\sum_{p\geq q-2} \Delta_q(\tilde u_p \cdot\nabla v_p).
\end{split}
\end{equation}
To reveal cancelations in the estimates, we also introduce the commutator
\begin{equation}\notag
[\Delta_q, u_{\leq{p-2}}\cdot\nabla]v_p=\Delta_q(u_{\leq{p-2}}\cdot\nabla v_p)-u_{\leq{p-2}}\cdot\nabla \Delta_qv_p.
\end{equation}
The following estimate holds.
\begin{Lemma} \label{le-comm1}
Let $u$ be a function with $\nabla\cdot u=0$. For any $1\leq r_1, r_2, r_3\leq \infty$ satisfying $\frac1{r_2}+\frac1{r_3}=\frac1{r_1}$, we have
\[\|[\Delta_q,u_{\leq{p-2}}\cdot\nabla] v_q\|_{r_1}\lesssim \|v_q\|_{r_3}\sum_{p' \leq p-2} \lambda_{p'} \|u_{p'}\|_{r_2}.\]
\end{Lemma}
\pf
Following the definition of $\Delta_q$,  we infer
\begin{equation}\notag
\begin{split}
[\Delta_q,u_{\leq{p-2}}\cdot\nabla] v_q=&\int_{\R^3}\lambda_q^3h(\lambda_q(x-y))\left(u_{\leq p-2}(y)-u_{\leq p-2}(x)\right)\cdot\nabla v_q(y)\,dy\\
=&-\int_{\R^3}\lambda_q^3\nabla h(\lambda_q(x-y))\left(u_{\leq p-2}(y)-u_{\leq p-2}(x)\right)\otimes v_q(y)\,dy,
\end{split}
\end{equation}
thanks to the fact $\nabla\cdot u_{\leq p-2}=0$.
Thus, by Young's inequality,
\begin{equation}\notag
\begin{split}
\|[\Delta_q,u_{\leq{p-2}}&\cdot\nabla] u_q\|_{r_1}\\
\lesssim &\|v_q\|_{r_3}\sum_{p' \leq p-2} \lambda_{p'} \|u_{p'}\|_{r_2}\left|\int_{\R^3}\lambda_q^3|x-y|\nabla h(\lambda_q(x-y))\, dy\right|\\
\lesssim &\|v_q\|_{r_3}\sum_{p' \leq p-2} \lambda_{p'} \|u_{p'}\|_{r_2},
\end{split}
\end{equation}
 for $\frac1{r_2}+\frac1{r_3}=\frac1{r_1}$.
 \cbdu
 
One can see another benefit of using commutator is that the derivative on high modes can be moved onto the low modes.
 We define a few more commutators regarding the $\Q$ tensor terms in the same spirit, as follows
\begin{equation}\notag
\begin{split}
[\Delta_q, \Delta\Q_p]\Q_{\leq p-2}=&\Delta_q(\Delta \Q_p\Q_{\leq p-2})-\Delta\Delta_q\Q_p\Q_{\leq p-2},\\
[\Delta_q, \Q_{\leq{p-2}}\Delta]\Q_p=&\Delta_q(\Q_{\leq p-2}\Delta \Q_p)-\Q_{\leq p-2}\Delta\Delta_q\Q_p,\\
[\Delta_q, \Omega(u)_p]\Q_{\leq p-2}=&\Delta_q(\Omega(u)_p\Q_{\leq p-2})-\Delta_q\Omega(u)_p\Q_{\leq p-2},\\
[\Delta_q, \Q_{\leq{p-2}}]\Omega(u)_p=&\Delta_q(\Q_{\leq p-2}\Omega(u)_p)-\Q_{\leq p-2}\Delta_q\Omega(u)_p,\\
[\Delta_q, \nabla\Q_p]\nabla\Q_{\leq p-2}=&\Delta_q(\nabla \Q_p\otimes\nabla\Q_{\leq p-2})-\nabla\Delta_q\Q_p\otimes\nabla\Q_{\leq p-2},\\
[\Delta_q, \nabla\Q_{\leq p-2}]\nabla\Q_p=&\Delta_q(\nabla \Q_{\leq p-2}\otimes\nabla\Q_p)-\nabla\Q_{\leq p-2}\otimes\nabla\Delta_q\Q_p.
\end{split}
\end{equation}
%We have the following estimates for these commutators.
\begin{Lemma}\label{le-comm2}
For any $1\leq r_1, r_2, r_3\leq \infty$ satisfying $\frac1{r_2}+\frac1{r_3}=\frac1{r_1}$, we have
\[\|[\Delta_q,  \Delta \Q_p]\Q_{\leq p-2}\|_{r_1}\lesssim \|\Q_p\|_{r_3}\sum_{p' \leq p-2} \lambda_{p'}^2 \|\Q_{p'}\|_{r_2},\]
\[\|[\Delta_q,  \Q_{\leq p-2}\Delta]\Q_p\|_{r_1}\lesssim \|\Q_p\|_{r_3}\sum_{p' \leq p-2} \lambda_{p'}^2 \|\Q_{p'}\|_{r_2},\]
\[\|[\Delta_q,  \Omega(u)_p]\Q_{\leq p-2}\|_{r_1}\lesssim \|u_p\|_{r_3}\sum_{p' \leq p-2} \lambda_{p'} \|\Q_{p'}\|_{r_2},\]
\[\|[\Delta_q, \Q_{\leq p-2}] \Omega(u)_p\|_{r_1}\lesssim \|u_p\|_{r_3}\sum_{p' \leq p-2} \lambda_{p'} \|\Q_{p'}\|_{r_2},\]
\[\|[\Delta_q,  \nabla \Q_p]\nabla\Q_{\leq p-2}\|_{r_1}\lesssim \|\Q_p\|_{r_3}\sum_{p' \leq p-2} \lambda_{p'}^2 \|\Q_{p'}\|_{r_2},\]
\[\|[\Delta_q,  \nabla \Q_{\leq p-2}]\nabla\Q_p\|_{r_1}\lesssim \|\Q_p\|_{r_3}\sum_{p' \leq p-2} \lambda_{p'}^2 \|\Q_{p'}\|_{r_2}.\]    
\end{Lemma}
\pf
Only the first inequality will be proven in the following, and other ones can be obtained in an analogous way.
Again, following the definition of $\Delta_q$ and applying integration by parts, the commutator can be written as
\begin{equation}\notag
\begin{split}
[\Delta_q,  \Delta \Q_p]\Q_{\leq p-2}=&\int_{\R^3}\lambda_q^3h(\lambda_q(x-y))\Delta\Q_p(y)\left(\Q_{\leq p-2}(y)-\Q_{\leq p-2}(x)\right)\,dy\\
%=&-\lambda_q^3\int_{\R^3}\nabla h(\lambda_q(x-y))\nabla\Q_p(y)\left(\Q_{\leq p-2}(y)-\Q_{\leq p-2}(x)\right) \,dy\\
%&-\lambda_q^3\int_{\R^3} h(\lambda_q(x-y))\nabla\Q_p(y)\nabla\left(\Q_{\leq p-2}(y)-\Q_{\leq p-2}(x)\right) \,dy\\
=&\int_{\R^3}\lambda_q^3\Delta h(\lambda_q(x-y))\Q_p(y)\left(\Q_{\leq p-2}(y)-\Q_{\leq p-2}(x)\right) \,dy\\
&+2\int_{\R^3}\lambda_q^3\nabla h(\lambda_q(x-y))\Q_p(y)\nabla\left(\Q_{\leq p-2}(y)-\Q_{\leq p-2}(x)\right) \,dy\\
%&+\int_{\R^3} \lambda_q^3\nabla h(\lambda_q(x-y))\Q_p(y)\nabla\left(\Q_{\leq p-2}(y)-\Q_{\leq p-2}(x)\right) \,dy\\
&+\int_{\R^3} \lambda_q^3h(\lambda_q(x-y))\Q_p(y)\Delta\left(\Q_{\leq p-2}(y)-\Q_{\leq p-2}(x)\right) \,dy.
\end{split}
\end{equation}
Thus, by Young's inequality, we infer
\begin{equation}\notag%\label{est-comm}
\begin{split}
&\|[\Delta_q,  \Delta \Q_p]\Q_{\leq p-2}\|_{r_1}\\
\lesssim &\|\Q_p\|_{r_3}\sum_{p' \leq p-2}  \|\nabla^2\Q_{p'}\|_{r_2}\left|\int_{\R^3}\lambda_q^3|x-y|^2\Delta h(\lambda_q(x-y))\, dy\right|\\
&+\|\Q_p\|_{r_3}\sum_{p' \leq p-2}\|\nabla^2\Q_{p'}\|_{r_2}\left|\int_{\R^3}\lambda_q^3|x-y|\nabla h(\lambda_q(x-y))\, dy\right|\\
&+\|\Q_p\|_{r_3}\sum_{p' \leq p-2}  \|\Delta\Q_{p'}\|_{r_2}\left|\int_{\R^3}\lambda_q^3h(\lambda_q(x-y))\, dy\right|\\
\lesssim &\|\Q_p\|_{r_3}\sum_{p' \leq p-2} \lambda_{p'}^2 \|\Q_{p'}\|_{r_2},
\end{split}
\end{equation}
 for $\frac1{r_2}+\frac1{r_3}=\frac1{r_1}$.  \cbdu

Similar computation strategy yields the following estimates. 
\begin{Lemma}\label{le-comm3}
For any $1\leq r_1, r_2, r_3\leq \infty$ satisfying $\frac1{r_2}+\frac1{r_3}=\frac1{r_1}$, we have
%\[\|\nabla([\Delta_q,  \Delta \Q_p]\Q_{\leq p-2})\|_{r_1}\lesssim \|\Q_p\|_{r_3}\sum_{p' \leq p-2} \lambda_{p'}^2 \|\Q_{p'}\|_{r_2},\]
% \[\|\nabla([\Delta_q,  \Q_{\leq p-2}\Delta]\Q_p)\|_{r_1}\lesssim \|\Q_p\|_{r_3}\sum_{p' \leq p-2} \lambda_{p'}^2 \|\Q_{p'}\|_{r_2},\]
\[\|\nabla([\Delta_q,  \Omega(u)_p]\Q_{\leq p-2})\|_{r_1}\lesssim \|u_p\|_{r_3}\sum_{p' \leq p-2} \lambda_{p'}^2 \|\Q_{p'}\|_{r_2},\]
\[\|\nabla([\Delta_q, \Q_{\leq p-2}] \Omega(u)_p)\|_{r_1}\lesssim \|u_p\|_{r_3}\sum_{p' \leq p-2} \lambda_{p'}^2 \|\Q_{p'}\|_{r_2}.\]
\end{Lemma}
%{\color{red} The last two can be used to estimate $J_{21}+J_{22}$, namely
%\[\|\nabla([\Delta_q, \nabla\Q_p] \nabla\Q_{\leq p-2})\|_{r_1}\lesssim \|\Q_p\|_{r_3}\sum_{p' \leq p-2} \lambda_{p'}^3 \|\Q_{p'}\|_{r_2},\]
%\[\|\nabla([\Delta_q, \nabla\Q_{\leq p-2}] \nabla\Q_p)\|_{r_1}\lesssim \|\Q_p\|_{r_3}\sum_{p' \leq p-2} \lambda_{p'}^3 \|\Q_{p'}\|_{r_2},\]}

\bigskip

\section{Regularity Criterion}
\label{sec:reg}

In this section we will establish the regularity criterion in Theorem \ref{thm}. Let $(u(t), \Q(t))$ be a weak solution of (\ref{LCD}) on $[0,T]$.
Based on the scaling of the Navier-Stokes equation, we define a dissipation wavenumber corresponding to velocity as
\begin{equation}\label{wave}
\Lambda(t)=\min \left\{\lambda_q:\lambda_p^{-1+\frac 3r}\|u_p(t)\|_r<c_r\min\{\nu,\mu\}, \forall p>q, q\in \mathbb{N} \right\},
\end{equation}
where $c_r$ is an adimensional  constant that depends only on $r$, and $r \in [2,6)$. We point out that the quantity $\lambda_p^{-1+\frac 3r}\|u_p(t)\|_r$ is scaling invariant. Let $Q(t)\in\mathbb N$ be such that $\lambda_{Q(t)}=\Lambda(t)$. It follows immediately that
\[
\|u_p(t)\|_r < \lambda_p^{1-\frac{3}{r}} c_r\min\{\nu,\mu\}, \qquad \forall p>Q(t),
\]
and 
\bg\label{Q}
\|u_{Q(t)}(t)\|_r\geq c_r\min\{\nu,\mu\}\Lambda^{1-\frac3r}(t),
\ed
provided $1<\Lambda(t)<\infty$. We also denote
\[
f(t)=\|\nabla u_{\leq Q(t)}\|_{B^0_{\infty,\infty}}.
\]
%\[f_{q^*}(t)=\|u_{(q^*, Q(t)]}\|_{B^1_{\infty,\infty}}=
%\left\{
%\begin{split}
%\sup_{\substack{q\leq Q(t)\\ q>q^*}}\lambda_q\|u_q(t)\|_\infty&, \qquad & q^* < Q(t),\\
%0&, \qquad &q^* \geq Q(t).
%\end{split}
%\right.
%\]

\subsection{Proof of Theorem \ref{thm}}
In order to prove that $(u,\Q)$ does not blow up at $T$, it is sufficient to show that $\|u(t)\|_{H^s}+\|\nabla\Q(t)\|_{H^s}$ is  bounded on $(0,T)$ for some $s>\frac{1}{2}$.  Multiplying equations  (\ref{LCD}) with $\Delta^2_qu$ and $\Delta^2_q\Delta\Q$ respectively yields 
%{\color{red} to have cancelations in $K$ and $J+L$, we may need to multiply the second equation by $\Delta_q^2\nabla^2 \Q$}
\begin{equation}\notag%\label{ineq:uq}
\begin{split}
\frac{1}{2}\frac{d}{dt}\|u_q\|_2^2\leq &-\nu\|\nabla u_q\|_2^2-\int_{\R^3}\Delta_q(u\cdot\nabla u)\cdot u_q\, dx
-\int_{\R^3}\Delta_q(\Sigma(\Q))\cdot \nabla u_q\, dx,\\
\frac{1}{2}\frac{d}{dt}\|\nabla \Q_q\|_2^2\leq &-\mu\|\Delta \Q_q\|_2^2-\int_{\R^3}\Delta_q(u\cdot\nabla \Q)\cdot \Delta\Q_qdx\\
&+\int_{\R^3}\Delta_q(\S(\nabla u, \Q))\cdot \Delta \Q_qdx-\int_{\R^3}\Delta_q\Ll[\partial F(\Q)] \Delta \Q_qdx.
%:=&-\mu\|\nabla b_q\|_2^2+I_3+I_4.
\end{split}
\end{equation}
Adding the above two inequalities, multiplying by $\lambda_q^{2s}$, and adding them for all $q\geq -1$, 
we obtain
\begin{equation}\label{ineq-energy}
\begin{split}
\frac{1}{2}\frac{d}{dt}\sum_{q\geq -1}\lambda_q^{2s}\left(\|u_q\|_2^2+\|\nabla \Q_q\|_2^2\right)\leq &-\sum_{q\geq -1}\lambda_q^{2s}\left(\nu\|\nabla u_q\|_2^2+\mu\|\Delta \Q_q\|_2^2\right)\\
&-(I+J+K+L+M),
\end{split}
\end{equation}
with
\begin{equation}\notag
\begin{split}
I=&\sum_{q\geq -1}\lambda_q^{2s}\int_{\R^3}\Delta_q(u\cdot\nabla u)\cdot u_q\, dx, \qquad
J=-\sum_{q\geq -1}\lambda_q^{2s}\int_{\R^3}\Delta_q(\Sigma(\Q)) \nabla u_q\, dx,\\
K=&\sum_{q\geq -1}\lambda_q^{2s}\int_{\R^3}\Delta_q(u\cdot\nabla \Q)\cdot \Delta\Q_q\, dx,\qquad
L=\sum_{q\geq -1}\lambda_q^{2s}\int_{\R^3}\Delta_q(\S(\nabla u,\Q))\cdot \Delta\Q_q\, dx,\\
M=&\sum_{q\geq -1}\lambda_q^{2s}\int_{\R^3}\Delta_q\Ll[\partial F(\Q)] \Delta \Q_q\, dx.
\end{split}
\end{equation}
Thanks to the maximum principle stated in Lemma \ref{max}, 
\[\|\Q(t,\cdot)\|_{L^\infty}\leq C \ \ \quad \mbox { for all } t\geq 0,\]
the term $M$ can be estimated immediately. 
Recalling
$F(\Q)=\frac a2|\Q|^2+\frac b3 \rm{tr}[\Q^3]+\frac c4|\Q|^4$, we have
\begin{equation}\label{est-m}
\begin{split}
|M|\leq &\sum_{q\geq -1}\lambda_q^{2s}\int_{\R^3}|\Delta_q\Ll[\partial F(\Q)] \Delta \Q_q|\, dx\\
\leq &\sum_{q\geq -1}\lambda_q^{2s}\int_{\R^3}|a\Q_q+b\Delta_q(\Q^2)+c\Delta_q(\Q{\rm{tr}}\Q^2)||\Delta \Q_q|\, dx\\
\lesssim&(1+\|\Q\|_\infty+\|\Q\|_\infty^2)\sum_{q\geq -1}\lambda_q^{2s+2}\|\Q_q\|_2^2\\
\lesssim&\sum_{q\geq -1}\lambda_q^{2s+2}\|\Q_q\|_2^2.
\end{split}
\end{equation}

Regarding the other terms, the main idea is to decompose them into high frequency and low frequency parts (by $Q$), such that the high frequency parts get absorbed by the diffusion term $\nu\|u\|_{H^{s+1}}^2+\mu\|\Q\|_{H^{s+2}}^2$. 
The term $I$ can be dealt with the same way as for the Navier-Stokes equation in \cite{CD}, and the estimate is
\begin{equation}\label{est-i}
|I|\lesssim c_r\mu\sum_{q>Q} \lambda_q^{2s+2}\|u_q\|_2^2+Qf(t)\sum_{q\geq -1} \lambda_q^{2s}\|u_q\|_2^2.
\end{equation}
%While the term $K$ can be estimated as the term $\int_{\R^3} (u\cdot\nabla) b\cdot b\, dx$ as for the MHD system in \cite{CD-mhd}, since $\nabla \Q$ has the same scaling as $b$. Also, there may be cancelation in $K+J_2$. Thus, 

We proceed the estimate for $J$, $L$ and $K$ in the following.
Recall 
\[\Sigma(\Q)=\Delta \Q\Q-\Q\Delta\Q-\nabla \Q\otimes\nabla\Q, \ \quad \S(\nabla u,\Q)=\Omega(u)\Q-\Q\Omega(u).\]
It follows that 
\begin{equation}\notag
\begin{split}
J=&-\sum_{q\geq -1}\lambda_q^{2s}\int_{\R^3}\Delta_q(\Delta \Q\Q-\Q\Delta\Q) \nabla u_q\, dx
+\sum_{q\geq -1}\lambda_q^{2s}\int_{\R^3}\Delta_q(\nabla \Q\otimes\nabla\Q) \nabla u_q\, dx\\
=&:J_1+J_2.
\end{split}
\end{equation}
We shall discover cancelations in $J_1+L$ and $J_2+K$ which are essential to obtain the ultimate estimate. Using Bony's paraproduct decomposition and the commutator notation, $J_1$ is decomposed as
\begin{equation}\notag%\label{eq:i2}
\begin{split}
J_1=
&\sum_{q\geq -1}\sum_{|q-p|\leq 2}\lambda_q^{2s}\int_{\R^3}\Delta_q(\Delta\Q_{p}\Q_{\leq{p-2}}-\Q_{\leq p-2}\Delta\Q_p)\nabla u_q\, dx\\
&+\sum_{q\geq -1}\sum_{|q-p|\leq 2}\lambda_q^{2s}\int_{\R^3}\Delta_q(\Delta\Q_{\leq p-2}\Q_p-\Q_p\Delta\Q_{\leq p-2})\nabla u_q\, dx\\
&+\sum_{q\geq -1}\sum_{p\geq q-2}\lambda_q^{2s}\int_{\R^3}\Delta_q(\Delta\Q_p\tilde \Q_p-\Q_p\Delta\tilde\Q_p)\nabla u_q\, dx\\
=&J_{11}+J_{12}+J_{13},
\end{split}
\end{equation}
with 
\begin{equation}\notag
\begin{split}
J_{11}=&\sum_{q\geq -1}\sum_{|q-p|\leq 2}\lambda_q^{2s}\int_{\R^3}\left([\Delta_q, \Delta\Q_p]\Q_{\leq p-2}-[\Delta_q,\Q_{\leq p-2}\Delta]\Q_p\right)\nabla u_q\, dx\\
&+\sum_{q\geq -1}\sum_{|q-p|\leq 2}\lambda_q^{2s}\int_{\R^3}\left(\Delta\Delta_q\Q_p\Q_{\leq q-2}-\Q_{\leq q-2}\Delta\Delta_q\Q_p\right)\nabla u_q\, dx\\
&+\sum_{q\geq -1}\sum_{|q-p|\leq 2}\lambda_q^{2s}\int_{\R^3}\left(\Delta\Delta_q\Q_p(\Q_{\leq{p-2}}-\Q_{\leq{q-2}})-(\Q_{\leq{p-2}}-\Q_{\leq{q-2}})\Delta\Delta_q\Q_p\right)\nabla u_q\, dx\\
=&J_{111}+J_{112}+J_{113}.
\end{split}
\end{equation}
%where we used {\color{red}$\sum_{|q-p|\leq 2}\Delta_pb_q=b_q$}.
%We will see that the term $J_{112}$ cancels a part of $L$. 
Similarly, $L$ can be decomposed as,
\begin{equation}\notag
\begin{split}
L
%=&-\sum_{q\geq -1}\lambda_q^{2s}\int_{\R^3}\Delta_q(\Omega(u)\Q-\Q\Omega(u))\cdot \Delta\Q_q\, dx\\
=&-\sum_{q\geq -1}\sum_{|q-p|\leq 2}\lambda_q^{2s}\int_{\R^3}\Delta_q(\Omega(u)_p\Q_{\leq p-2}-\Q_{\leq p-2}\Omega(u)_p)\cdot \Delta\Q_q\, dx\\
&-\sum_{q\geq -1}\sum_{|q-p|\leq 2}\lambda_q^{2s}\int_{\R^3}\Delta_q(\Omega(u)_{\leq p-2}\Q_p-\Q_p\Omega(u)_{\leq p-2})\cdot \Delta\Q_q\, dx\\
&-\sum_{q\geq -1}\sum_{p\geq q-2}\lambda_q^{2s}\int_{\R^3}\Delta_q(\Omega(u)_p\tilde\Q_p-\Q_p\Omega(\tilde u)_p)\cdot \Delta\Q_q\, dx\\
=&L_1+L_2+L_3,
\end{split}
\end{equation}
with 
\begin{equation}\notag
\begin{split}
L_1=&-\sum_{q\geq -1}\sum_{|q-p|\leq 2}\lambda_q^{2s}\int_{\R^3}\left([\Delta_q, \Omega(u)_p]\Q_{\leq p-2}-[\Delta_q, \Q_{\leq p-2}]\Omega(u)_p\right)\cdot \Delta\Q_q\, dx\\
&-\sum_{q\geq -1}\sum_{|q-p|\leq 2}\lambda_q^{2s}\int_{\R^3}\left(\Delta_q \Omega(u)_p\Q_{\leq q-2}-\Q_{\leq q-2}\Delta_q\Omega(u)_p\right)\cdot \Delta\Q_q\, dx\\
&-\sum_{q\geq -1}\sum_{|q-p|\leq 2}\lambda_q^{2s}\int_{\R^3}\left(\Delta_q\Omega(u)_p(\Q_{\leq p-2}-\Q_{q-2})-(\Q_{\leq p-2}-\Q_{\leq q-2})\Delta_q\Omega(u)_p\right)\cdot \Delta\Q_q\, dx\\
=&L_{11}+L_{12}+L_{13}.
\end{split}
\end{equation}
Note that $\sum_{|q-p|\leq 2}\Delta\Delta_q\Q_p=\Delta\Q_q$ and $\sum_{|q-p|\leq 2}\Delta_q \Omega(u)_p=\Omega(u)_q$. Thus 
\begin{equation}\notag
\begin{split}
J_{112}+L_{12}=&\sum_{q\geq -1}\lambda_q^{2s}\int_{\R^3}\left(\Delta\Q_q\Q_{\leq q-2}-\Q_{\leq q-2}\Delta\Q_q\right)\nabla u_q\, dx\\
&-\sum_{q\geq -1}\lambda_q^{2s}\int_{\R^3}\left(\Omega(u)_q\Q_{\leq q-2}-\Q_{\leq q-2}\Omega(u)_q\right)\cdot \Delta\Q_q\, dx\\
=&0.
\end{split}
\end{equation}
The other terms in $J_1$ and $L$ are estimated as follows.
We further split $J_{111}$ into high and low frequency parts as 
\begin{equation}\notag
\begin{split}
|J_{111}|\leq &2 \sum_{q>Q}\sum_{|q-p|\leq 2}\lambda_q^{2s}\int_{\R^3}\left|[\Delta_q, \Delta\Q_p] \Q_{\leq p-2}\nabla u_q\right|\, dx\\
&+2\sum_{-1\leq q\leq Q}\sum_{|q-p|\leq 2}\lambda_q^{2s}\int_{\R^3}\left|[\Delta_q, \Delta\Q_p] \Q_{\leq p-2}\nabla u_q\right|\, dx\\
\equiv & J_{1111}+J_{1112}.
\end{split}
\end{equation}
By H\"older's inequality, the commutator estimate in Lemma \ref{le-comm2},  the definition of $\Lambda_r$, and Jensen's inequality, it follows that  
\begin{equation}\notag
\begin{split}
J_{1111}
\lesssim&\sum_{q> Q}\lambda_q^{2s}\|\nabla u_q\|_r\sum_{|q-p|\leq 2}\| \Q_p\|_2\sum_{p'\leq p-2}\lambda_{p'}^2\|\Q_{p'}\|_{\frac{2r}{r-2}}\\
\lesssim&c_r\mu \sum_{q> Q}\lambda_q^{2s+2-\frac3r}\sum_{|q-p|\leq 2}\| \Q_p\|_2\sum_{p'\leq q}\lambda_{p'}^{2+\frac3r}\|\Q_{p'}\|_2\\
\lesssim&c_r\mu \sum_{q> Q-2}\lambda_q^{2s+2-\frac3r}\|\Q_q\|_2\sum_{p'\leq q}\lambda_{p'}^{2+\frac3r}\|\Q_{p'}\|_2\\
\lesssim&c_r\mu\sum_{q> Q-2}\lambda_q^{s+2}\|\Q_q\|_2\sum_{p'\leq q}\lambda_{p'}^{s+2}\|\Q_{p'}\|_2\lambda_{q-p'}^{s-\frac 3r}\\
\lesssim&c_r\mu\sum_{q\geq -1}\lambda_q^{2s+4}\|\Q_q\|_2^2,
\end{split}
\end{equation}
where we needed  $2\leq r<\frac 3s$ 
%{\color{red}Adding a kernel to the definition of $\Lambda$ may relax the range for $r$}.
Also, by H\"older's inequality, the definition of $f(t)$ and Jensen's inequality, 
\begin{equation}\notag
\begin{split}
J_{1112}
\lesssim &\sum_{-1\leq q\leq Q}\lambda_q^{2s}\|\nabla u_q\|_\infty\sum_{|q-p|\leq 2}\|\Q_p\|_2\sum_{p'\leq p-2}\lambda_{p'}^2\|\Q_{p'}\|_2\\
\lesssim &f(t)\sum_{-1\leq q\leq Q}\lambda_q^{2s}\sum_{|q-p|\leq 2}\|\Q_p\|_2\sum_{p'\leq q}\lambda_{p'}^2\|\Q_{p'}\|_2\\
\lesssim&f(t)\sum_{-1\leq q\leq Q+2}\lambda_q^{2s}\|\Q_q\|_2\sum_{p'\leq q}\lambda_{p'}^2\|\Q_{p'}\|_2\\
\lesssim&f(t)\sum_{-1\leq q\leq Q+2}\lambda_q^{s+1}\|\Q_q\|_2\sum_{p'\leq q}\lambda_{p'}^{s+1}\|\Q_{p'}\|_2\lambda_{p'-q}^{1-s}\\
\lesssim&f(t)\sum_{-1\leq q\leq Q+2}\lambda_q^{2s+2}\|\Q_q\|_2^2,
\end{split}
\end{equation}
where we used $\frac 12< s<1$.  
Similar analysis yields 
\begin{equation}\notag
\begin{split}
|J_{113}|\leq & \sum_{q\geq -1}\sum_{|q-p|\leq 2}\lambda_q^{2s}\int_{\R^3}\left(\Delta\Delta_q\Q_p(\Q_{\leq{p-2}}-\Q_{\leq{q-2}})-(\Q_{\leq{p-2}}-\Q_{\leq{q-2}})\Delta\Delta_q\Q_p\right)\nabla u_q\, dx\\
\leq & \sum_{q>Q}\sum_{|q-p|\leq 2}\lambda_q^{2s+3}\|u_q\|_r\|\Q_{\leq{p-2}}-\Q_{\leq{q-2}}\|_2\|\Q_p\|_{\frac{2r}{r-2}}\\
& +\sum_{-1\leq q\leq Q}\sum_{|q-p|\leq 2}\lambda_q^{2s+3}\|u_q\|_\infty\|\Q_{\leq{p-2}}-\Q_{\leq{q-2}}\|_2\|\Q_p\|_2\\
\lesssim & c_r\mu \sum_{q>Q}\lambda_q^{2s+4-\frac3r}\sum_{|q-p|\leq 2}\lambda_p^{\frac3r}\|\Q_{\leq{p-2}}-\Q_{\leq{q-2}}\|_2\|\Q_p\|_2\\
& +f(t)\sum_{-1\leq q\leq Q}\lambda_q^{2s+2}\sum_{|q-p|\leq 2}\|\Q_{\leq{p-2}}-\Q_{\leq{q-2}}\|_2\|\Q_p\|_2\\
\lesssim & c_r\mu \sum_{q>Q}\lambda_q^{2s+4-\frac3r}\sum_{q-3\leq p\leq q+2}\lambda_p^{\frac3r}\|\Q_p\|_2^2 +f(t)\sum_{-1\leq q\leq Q}\lambda_q^{2s+2}\|\Q_q\|_2^2\\
\lesssim & c_r\mu \sum_{q>Q-3}\lambda_q^{2s+4}\|\Q_p\|_2^2 +f(t)\sum_{-1\leq q\leq Q}\lambda_q^{2s+2}\|\Q_q\|_2^2.
\end{split}
\end{equation}
Notice that $J_{12}$ enjoys the same estimate as $J_{111}$. While for $J_{13}$ we have
\begin{equation}\notag
\begin{split}
|J_{13}|\lesssim &\sum_{q\geq -1}\sum_{p\geq q-2}\lambda_q^{2s}\int_{\R^3}\left|\Delta_q(\Delta\Q_p\tilde \Q_p-\Q_p\Delta\tilde\Q_p)\nabla u_q\right|\, dx\\
%\lesssim &\sum_{q\geq -1}\sum_{p\geq q-2}\lambda_q^{2s}\int_{\mathbb R^3}|\Delta_q(b_p\otimes\tilde b_p)\nabla u_q|\, dx\\
\lesssim &\sum_{q> Q}\lambda_q^{2s+1}\|u_q\|_\infty\sum_{p\geq q-2}\lambda_p^2\|\Q_p\|_2^2
+\sum_{-1\leq q\leq Q}\lambda_q^{2s+1}\|u_q\|_\infty\sum_{p\geq q-2}\lambda_p^2\|\Q_p\|_2^2\\
\lesssim &\sum_{q> Q}\lambda_q^{2s+1+\frac3r}\|u_q\|_r\sum_{p\geq q-2}\lambda_p^2\|\Q_p\|_2^2
+f(t)\sum_{-1\leq q\leq Q}\lambda_q^{2s}\sum_{p\geq q-2}\lambda_p^2\|\Q_p\|_2^2\\
\lesssim &c_r\mu\sum_{p> Q}\lambda_p^{2s+4}\|\Q_p\|_2^2\sum_{Q< q\leq p+2}\lambda_{q-p}^{2s+2}
+f(t)\sum_{p\geq -1}\lambda_p^{2s+2}\|\Q_p\|_2^2\sum_{-1\leq q\leq p}\lambda_{q-p}^{2s}\\
\lesssim &c_r\mu\sum_{q> Q}\lambda_q^{2s+4}\|\Q_q\|_2^2+f(t)\sum_{q\geq -1}\lambda_q^{2s+2}\|\Q_q\|_2^2.
\end{split}
\end{equation}
Therefore, for $2\leq r<\frac 3s$ and $\frac 12<s<1$, 
\begin{equation}\label{est-j1}
|J_1|\lesssim c_r\mu\sum_{q> Q-3}\lambda_q^{2s+4}\|\Q_q\|_2^2+f(t)\sum_{q\geq -1}\lambda_q^{2s+2}\|\Q_q\|_2^2.
\end{equation}
%Now we estimate the other terms in $L$.
Applying integration by parts to $L_{11}$ yields
\begin{equation}\notag
\begin{split}
L_{11}= &-\sum_{p\leq Q}\sum_{|q-p|\leq 2}\lambda_q^{2s}\int_{\R^3}\nabla\left([\Delta_q, \Omega(u)_p]\Q_{\leq p-2}-[\Delta_q, \Q_{\leq p-2}]\Omega(u)_p\right)\nabla\Q_q\, dx\\
&-\sum_{p>Q}\sum_{|q-p|\leq 2}\lambda_q^{2s}\int_{\R^3}\nabla\left([\Delta_q, \Omega(u)_p]\Q_{\leq p-2}-[\Delta_q, \Q_{\leq p-2}]\Omega(u)_p\right)\nabla\Q_q\, dx\\
\equiv& L_{111}+L_{112}.
\end{split}
\end{equation}
By Lemma \ref{le-comm3}, Young's inequality and Jensen's inequality, we infer
\begin{equation}\notag
\begin{split}
|L_{111}|\lesssim &\sum_{p\leq Q}\sum_{|q-p|\leq 2}\lambda_p^{2s}\|u_p\|_\infty\| \nabla\Q_q\|_2\sum_{p'\leq p-2}\lambda_{p'}^2\|\Q_{p'}\|_2\\
\lesssim &f(t)\sum_{p\leq Q}\lambda_p^{2s}\| \Q_p\|_2\sum_{p'\leq p-2}\lambda_{p'}^{2}\| \Q_{p'}\|_2\\
\lesssim &f(t)\sum_{p\leq Q}\lambda_p^{s+1}\| \Q_p\|_2\sum_{p'\leq p-2}\lambda_{p'}^{s+1}\| \Q_{p'}\|_2\lambda_{p-p'}^{s-1}\\
\lesssim &f(t)\sum_{p\leq Q}\lambda_p^{2s+2}\| \Q_p\|_2^2
\end{split}
\end{equation}
for $s<1$; while
\begin{equation}\notag
\begin{split}
|L_{112}|\lesssim &\sum_{p> Q}\sum_{|q-p|\leq 2}\lambda_p^{2s}\|u_p\|_r\| \nabla\Q_q\|_2\sum_{p'\leq p-2}\lambda_{p'}^2\|\Q_{p'}\|_{\frac{2r}{r-2}}\\
\lesssim &c_r\mu\sum_{p> Q}\lambda_p^{2s+2-\frac3r}\| \Q_p\|_2\sum_{p'\leq p-2}\lambda_{p'}^{2+\frac3r}\| \Q_{p'}\|_2\\
\lesssim &c_r\mu\sum_{p> Q}\lambda_p^{s+2}\| \Q_p\|_2\sum_{p'\leq p-2}\lambda_{p'}^{s+2}\| \Q_{p'}\|_2\lambda_{p'-p}^{\frac3r-s}\\
\lesssim &c_r\mu\sum_{p\geq -1}\lambda_p^{2s+4}\| \Q_p\|_2^2
\end{split}
\end{equation}
for $s<\frac3r$. 

Notice that $L_{13}$ can be estimated similarly to $J_{113}$ and $L_{3}$ can be estimated as $J_{13}$. Thus
\begin{equation}\notag
|L_{13}|+|L_{3}|\lesssim c_r\mu\sum_{q\geq -1}\lambda_q^{2s+4}\|\Q_q\|_2^2+f(t)\sum_{-1\leq q\leq Q}\lambda_q^{2s+2}\|\Q_q\|_2^2.
\end{equation}
To estimate $L_{2}$, we split the summation first as
\begin{equation}\notag
\begin{split}
|L_{2}|\leq&\sum_{q\geq -1}\sum_{|q-p|\leq 2}\lambda_q^{2s}\int_{\R^3}\left|\Delta_q(\Omega(u)_{\leq p-2}\Q_p-\Q_p\Omega(u)_{\leq p-2})\cdot \Delta\Q_q\right|\, dx\\
\lesssim &\sum_{-1\leq p\leq Q+2}\sum_{|p-q|\leq 2}\lambda_q^{2s}\int_{\mathbb R^3}\left|\Delta_q(\Omega(u)_{\leq p-2}\Q_p)\Delta\Q_q\right|\,dx\\
&+\sum_{p>Q+2}\sum_{|p-q|\leq 2}\lambda_q^{2s}\int_{\mathbb R^3}\left|\Delta_q(\Omega(u)_{\leq Q}\Q_p)\Delta\Q_q\right|\,dx\\
&+\sum_{p>Q+2}\sum_{|p-q|\leq 2}\lambda_q^{2s}\int_{\mathbb R^3}\left|\Delta_q(\Omega(u)_{(Q,p-2]}\Q_p)\Delta\Q_q\right|\,dx\\
\equiv & L_{21}+L_{22}+L_{23}.
\end{split}
\end{equation}
Then using H\"older's inequality and the definition of $f(t)$ we obtain
\begin{equation}\notag
\begin{split}
L_{21}\lesssim & \sum_{-1\leq p\leq Q+2}\|\nabla u_{\leq p-2}\|_\infty\|\Q_p\|_2\sum_{|p-q|\leq 2}\lambda_q^{2s}\|\Delta\Q_q\|_2\\
\lesssim & Qf(t)\sum_{-1\leq p\leq Q+2}\|\Q_p\|_2\sum_{|p-q|\leq 2}\lambda_q^{2s+2}\|\Q_q\|_2\\
\lesssim & Qf(t)\sum_{q\geq -1}\lambda_q^{2s+2}\|\Q_q\|_2^2,
\end{split}
\end{equation}
and
\begin{equation}\notag
\begin{split}
L_{22} \lesssim &\sum_{p> Q+2}\|\nabla u_{\leq Q}\|_\infty\|\Q_p\|_2\sum_{|p-q|\leq 2}\lambda_q^{2s}\|\Delta\Q_q\|_2\\
\lesssim & Qf(t)\sum_{p> Q+2}\|\Q_p\|_2\sum_{|p-q|\leq 2}\lambda_q^{2s+2}\|\Q_q\|_2\\
\lesssim & Qf(t)\sum_{q> Q}\lambda_q^{2s+2}\|\Q_q\|_2^2.
\end{split}
\end{equation}
Applying H\"older's inequality, the definition of $\Lambda_r$ and Jensen's inequality yields
\begin{equation}\notag
\begin{split}
L_{23} \lesssim &\sum_{p> Q+2}\|\nabla u_{(Q,p-2]}\|_r\|\Q_p\|_{\frac{2r}{r-2}}\sum_{|p-q|\leq 2}\lambda_q^{2s}\|\Delta\Q_q\|_2\\
\lesssim &\sum_{p> Q+2}\lambda_p^{\frac3r}\|\Q_p\|_2\sum_{|p-q|\leq 2}\lambda_q^{2s+2}\|\Q_q\|_2\sum_{Q<p'\leq p-2}\lambda_{p'}\|u_{p'}\|_r\\
\lesssim & c_r\mu \sum_{p> Q}\lambda_p^{2s+2+\frac3r}\|\Q_p\|_2^2\sum_{Q<p'\leq p-2}\lambda_{p'}^{2-\frac3r}\\
\lesssim & c_r\mu \sum_{p> Q}\lambda_p^{2s+4}\|\Q_p\|_2^2\sum_{Q<p'\leq p-2}\lambda_{p'-p}^{2-\frac3r}\\
\lesssim & c_r\mu \sum_{p> Q}\lambda_p^{2s+4}\|\Q_p\|_2^2,
\end{split}
\end{equation}
since $r\geq 2$. Thus we have established that 
\begin{equation}\label{est-l}
|L|\lesssim c_r\mu\sum_{q\geq -1}\lambda_q^{2s+4}\|\Q_q\|_2^2+Qf(t)\sum_{q\geq -1}\lambda_q^{2s+2}\|\Q_q\|_2^2.
\end{equation}

We deal with $J_2+K$ now. As before, $J_2$ can be decomposed as 
\begin{equation}\notag
\begin{split}
J_2=&\sum_{q\geq -1}\lambda_q^{2s}\int_{\R^3}\Delta_q(\nabla \Q\otimes\nabla\Q) \nabla u_q\, dx\\
=&\sum_{q\geq -1}\sum_{|q-p|\leq 2}\lambda_q^{2s}\int_{\R^3}\Delta_q(\nabla\Q_{\leq p-2}\otimes\nabla \Q_p)\nabla u_q\, dx\\
&+\sum_{q\geq -1}\sum_{|q-p|\leq 2}\lambda_q^{2s}\int_{\R^3}\Delta_q(\nabla\Q_{p}\otimes\nabla \Q_{\leq{p-2}})\nabla u_q\, dx\\
&+\sum_{q\geq -1}\sum_{p\geq q-2}\lambda_q^{2s}\int_{\R^3}\Delta_q(\nabla\Q_p\otimes\nabla\tilde \Q_p)\nabla u_q\, dx\\
=&J_{21}+J_{22}+J_{23}.
\end{split}
\end{equation}
Applying the commutator notation, we can rewrite $J_{21}$ and $J_{22}$
\begin{equation}\notag
\begin{split}
J_{21}=&\sum_{q\geq -1}\sum_{|q-p|\leq 2}\lambda_q^{2s}\int_{\R^3}[\Delta_q,\nabla\Q_{\leq p-2}]\nabla\Q_p\nabla u_q\, dx\\
&+\sum_{q\geq -1}\sum_{|q-p|\leq 2}\lambda_q^{2s}\int_{\R^3}\nabla\Q_{\leq q-2}\otimes\nabla\Delta_q\Q_p\nabla u_q\, dx\\
&+\sum_{q\geq -1}\sum_{|q-p|\leq 2}\lambda_q^{2s}\int_{\R^3}(\nabla\Q_{\leq p-2}-\nabla\Q_{\leq q-2})\otimes\nabla\Delta_q\Q_p\nabla u_q\, dx\\
\equiv&J_{211}+J_{212}+J_{213};
\end{split}
\end{equation}
and
\begin{equation}\notag
\begin{split}
J_{22}=&\sum_{q\geq -1}\sum_{|q-p|\leq 2}\lambda_q^{2s}\int_{\R^3}[\Delta_q,\nabla\Q_p]\nabla\Q_{\leq p-2}\nabla u_q\, dx\\
&+\sum_{q\geq -1}\sum_{|q-p|\leq 2}\lambda_q^{2s}\int_{\R^3}\nabla\Delta_q\Q_p\otimes\nabla\Q_{\leq q-2}\nabla u_q\, dx\\
&+\sum_{q\geq -1}\sum_{|q-p|\leq 2}\lambda_q^{2s}\int_{\R^3}\nabla\Delta_q\Q_p\otimes(\nabla\Q_{\leq p-2}-\nabla\Q_{\leq q-2})\nabla u_q\, dx\\
\equiv&J_{221}+J_{222}+J_{223}.
\end{split}
\end{equation}
Since $\sum_{|p-q|\leq 2}\nabla\Delta_q\Q_p=\nabla \Q_q$, we have
\begin{equation}\notag
J_{212}+J_{222}=\sum_{q\geq -1}\lambda_q^{2s}\int_{\R^3}\nabla\Q_{\leq q-2}\otimes\nabla\Q_q\nabla u_q\, dx
+\sum_{q\geq -1}\lambda_q^{2s}\int_{\R^3}\nabla\Q_q\otimes\nabla\Q_{\leq q-2}\nabla u_q\, dx
\end{equation}
which will be estimated together with $K_{22}$ later.

We also decompose $K$ by Bony's paraproduct,
\begin{equation}\notag%\label{eq:i3}
\begin{split}
K=&\sum_{q\geq -1}\lambda_q^{2s}\int_{\R^3}\Delta_q(u\cdot\nabla \Q)\cdot \Delta\Q_q\, dx\\
=&\sum_{q\geq -1}\sum_{|q-p|\leq 2}\lambda_q^{2s}\int_{\R^3}\Delta_q(u_{\leq p-2}\cdot\nabla \Q_p)\cdot\Delta\Q_q\, dx\\
&+\sum_{q\geq -1}\sum_{|q-p|\leq 2}\lambda_q^{2s}\int_{\R^3}\Delta_q(u_{p}\cdot\nabla \Q_{\leq{p-2}})\cdot\Delta\Q_q\, dx\\
&+\sum_{q\geq -1}\sum_{p\geq q-2}\lambda_q^{2s}\int_{\R^3}\Delta_q(u_p\cdot\nabla\tilde \Q_p)\cdot\Delta\Q_q\, dx\\
=&K_{1}+K_{2}+K_{3},
\end{split}
\end{equation}
with 
\begin{equation}\notag
\begin{split}
K_{1}=&\sum_{q\geq -1}\sum_{|q-p|\leq 2}\lambda_q^{2s}\int_{\R^3}[\Delta_q,u_{\leq{p-2}}\cdot\nabla] \Q_p\cdot\Delta\Q_q\, dx\\
&+\sum_{q\geq -1}\lambda_q^{2s}\int_{\R^3}(u_{\leq{q-2}}\cdot\nabla) \Q_q\cdot\Delta\Q_q\, dx\\
&+\sum_{q\geq -1}\sum_{|q-p|\leq 2}\lambda_q^{2s}\int_{\R^3}(u_{\leq{p-2}}-u_{\leq{q-2}})\cdot\nabla\Delta_q\Q_p\cdot\Delta\Q_q\, dx\\
=&K_{11}+K_{12}+K_{13};
\end{split}
\end{equation}
and
\begin{equation}\notag
\begin{split}
K_2=&\sum_{q\geq -1}\sum_{|q-p|\leq 2}\lambda_q^{2s}\int_{\R^3}\Delta_q(u_{p}\cdot\nabla \Q_{\leq{p-2}})\cdot\Delta\Q_q\, dx\\
=&\sum_{q\geq -1}\sum_{|q-p|\leq 2}\lambda_q^{2s}\int_{\R^3}[\Delta_q, u_{p}\cdot\nabla] \Q_{\leq{p-2}}\cdot\Delta\Q_q\, dx\\
&+\sum_{q\geq -1}\lambda_q^{2s}\int_{\R^3}(u_q\cdot\nabla) \Q_{\leq{q-2}}\cdot\Delta\Q_q\, dx\\
&+\sum_{q\geq -1}\sum_{|q-p|\leq 2}\lambda_q^{2s}\int_{\R^3}(\Delta_qu_{p}\cdot\nabla) (\Q_{\leq{p-2}}-\Q_{\leq{q-2}})\cdot\Delta\Q_q\, dx\\
\equiv &K_{21}+K_{22}+K_{23}.
\end{split}
\end{equation}
Here we used $\sum_{|q-p|\leq 2}\nabla\Delta_p\Q_q=\Q_q$ and $\sum_{|p-q|\leq 2}\Delta_qu_p=u_q$ to obtain $K_{12}$ and $K_{22}$, respectively. 

We claim that 
\begin{equation}\label{cancel3}
J_{212}+J_{222}+K_{22}=-\sum_{q\geq-1}\lambda_q^{2s}\int_{\R^3}(u_q\cdot\nabla)\Q_q\Delta\Q_{\leq q-2}\,dx.
\end{equation}
Indeed, denote $\A=\Q_{\leq q-2}$ and $\B=\Q_q$. Applying integration by parts and $\nabla\cdot u_q=0$ yields
\begin{equation}\notag
\begin{split}
&\int_{\R^3}(\nabla\Q_{\leq q-2}\otimes\nabla\Q_q)\nabla u_q\,dx+\int_{\R^3}(\nabla\Q_q\otimes\nabla\Q_{\leq q-2})\nabla u_q\,dx\\
=&\int_{\R^3}\A_{\gamma\delta,\alpha}\B_{\gamma\delta,\beta}(u_q)_{\alpha,\beta}\,dx
+\int_{\R^3}\B_{\gamma\delta,\alpha}\A_{\gamma\delta,\beta}(u_q)_{\alpha,\beta}\,dx\\
=&-\int_{\R^3}\A_{\gamma\delta,\alpha\beta}\B_{\gamma\delta,\beta}(u_q)_{\alpha}\,dx
-\int_{\R^3}\A_{\gamma\delta,\alpha}\B_{\gamma\delta,\beta\beta}(u_q)_{\alpha}\,dx\\
&-\int_{\R^3}\B_{\gamma\delta,\alpha\beta}\A_{\gamma\delta,\beta}(u_q)_{\alpha}\,dx
-\int_{\R^3}\B_{\gamma\delta,\alpha}\A_{\gamma\delta,\beta\beta}(u_q)_{\alpha}\,dx\\
=&-\int_{\R^3}\A_{\gamma\delta,\alpha\beta}\B_{\gamma\delta,\beta}(u_q)_{\alpha}\,dx
-\int_{\R^3}\A_{\gamma\delta,\alpha}\B_{\gamma\delta,\beta\beta}(u_q)_{\alpha}\,dx\\
&+\int_{\R^3}\B_{\gamma\delta,\beta}\A_{\gamma\delta,\alpha\beta}(u_q)_{\alpha}\,dx
+\int_{\R^3}\B_{\gamma\delta,\beta}\A_{\gamma\delta,\beta}(u_q)_{\alpha,\alpha}\,dx
-\int_{\R^3}\B_{\gamma\delta,\alpha}\A_{\gamma\delta,\beta\beta}(u_q)_{\alpha}\,dx\\
=&-\int_{\R^3}(u_q\cdot\nabla) \Q_{\leq{q-2}}\cdot\Delta\Q_q\, dx-
\int_{\R^3}(u_q\cdot\nabla)\Q_q\Delta\Q_{\leq q-2}\,dx.
\end{split}
\end{equation}
Thus, identity (\ref{cancel3}) follows immediately. 

We start now estimating the terms in $J_2+K$. It follows from (\ref{cancel3}) that
\begin{equation}\notag
\begin{split}
&|J_{212}+J_{222}+K_{22}|\\
\leq &\sum_{q>Q}\lambda_q^{2s}\int_{\R^3}|(u_q\cdot\nabla)\Q_q\Delta\Q_{\leq q-2}|\,dx+\sum_{q\leq Q}\lambda_q^{2s}\int_{\R^3}|(u_q\cdot\nabla)\Q_q\Delta\Q_{\leq q-2}|\,dx\\
\lesssim&\sum_{q>Q}\lambda_q^{2s}\|u_q\|_r\|\nabla \Q_q\|_2\sum_{p\leq q-2}\|\Delta\Q_p\|_{\frac{2r}{r-2}}+\sum_{q\leq Q}\lambda_q^{2s}\|u_q\|_\infty\|\nabla \Q_q\|_2\sum_{p\leq q-2}\|\Delta\Q_p\|_2\\
\lesssim&c_r\mu\sum_{q>Q}\lambda_q^{2s+2-\frac3r}\|\Q_q\|_2\sum_{p\leq q-2}\lambda_p^{2+\frac3r}\|\Q_p\|_2+f(t)\sum_{q\leq Q}\lambda_q^{2s}\|\Q_q\|_2\sum_{p\leq q-2}\lambda_p^2\|\Q_p\|_2\\
\lesssim&c_r\mu\sum_{q>Q}\lambda_q^{s+2}\|\Q_q\|_2\sum_{p\leq q-2}\lambda_p^{s+2}\|\Q_p\|_2\lambda_{q-p}^{s-\frac3r}+f(t)\sum_{q\leq Q}\lambda_q^{s+1}\|\Q_q\|_2\sum_{p\leq q-2}\lambda_p^{s+1}\|\Q_p\|_2\lambda_{q-p}^{s-1}\\
\lesssim&c_r\mu\sum_{q>Q}\lambda_q^{2s+4}\|\Q_q\|_2^2+f(t)\sum_{q\leq Q}\lambda_q^{2s+2}\|\Q_q\|_2^2
\end{split}
\end{equation}
for $s<\frac3r$ and $s<1$.

Some of the rest terms are relatively easy. For instance, $J_{211}+J_{221}$ can be estimated similarly as $J_{111}$, while $J_{213}+J_{223}$ similarly as $J_{113}$, $K_{13}$ and $K_{23}$ similarly as $L_{13}$. Also, recalling that one of the benefits of commutator is to move derivatives onto low frequency parts (see Lemma \ref{le-comm2}), we observe $K_{11}$ and $K_{21}$ can be handled in an analogous way as $L_2$ and $L_{11}$, respectively.

The term $J_{23}$ is estimated as follows,  
\begin{equation}\notag
\begin{split}
|J_{23}|\leq&\left|\sum_{q\geq -1}\lambda_q^{2s}\sum_{p\geq q-2}\int_{\mathbb R^3}\Delta_q(\nabla\Q_p\otimes\nabla\tilde\Q_p) \nabla u_q \, dx\right|\\
\leq&\sum_{q>Q}\lambda_q^{2s}\sum_{p\geq q-2}\int_{\mathbb R^3}\left|\Delta_q(\nabla\Q_p\otimes\nabla\tilde\Q_p) \nabla u_q\right| \, dx\\
&+\sum_{-1\leq q\leq Q}\lambda_q^{2s}\sum_{p\geq q-2}\int_{\mathbb R^3}\left|\Delta_q(\nabla\Q_p\otimes\nabla\tilde\Q_p) \nabla u_q\right| \, dx\\
\equiv & J_{231}+J_{232},
\end{split}
\end{equation}
with
\begin{equation}\notag
\begin{split}
|J_{231}|
\lesssim &\sum_{q> Q}\|\nabla u_q\|_\infty\sum_{p\geq q-2}\lambda_p^{2s}\|\nabla \Q_p\|_2^2\\
\lesssim &c_r\mu\sum_{q> Q}\lambda_q^2\sum_{p\geq q-2}\lambda_p^{2s+2}\|\Q_p\|_2^2\\
\lesssim &c_r\mu\sum_{p>Q-2}\lambda_p^{2s+4}\|\Q_p\|_2^2\sum_{q\leq p+2}\lambda_{q-p}^2\\
\lesssim &c_r\mu\sum_{q\geq -1}\lambda_q^{2s+4}\|\Q_q\|_2^2;
\end{split}
\end{equation}
and 
\begin{equation}\notag
\begin{split}
|J_{232}|
\lesssim &\sum_{-1\leq q\leq Q}\|\nabla u_q\|_\infty\sum_{p\geq q-2}\lambda_p^{2s}\|\nabla\Q_p\|_2^2\\
\lesssim &f(t)\sum_{-1\leq q\leq Q}\sum_{p\geq q-2}\lambda_p^{2s+2}\|\Q_p\|_2^2\\
\lesssim &f(t)\sum_{q\geq -1}\lambda_q^{2s+2}\|\Q_q\|_2^2.
\end{split}
\end{equation}

Noticing the cancelation  
\begin{equation}\notag%\label{cancel}
\int_{\R^3}(u\cdot\nabla )\Q\cdot\Delta\Q\,dx+\int_{\R^3}(\nabla\Q\otimes\nabla\Q)\nabla u\,dx=0, 
\end{equation}
%see the proof of Proposition 1 in \cite{PZ}. which will be utilized to estimate the most difficult part from $K$.
one can write $K_{12}=-\int_{\R^3}(\nabla\Q_q\otimes\nabla\Q_q)\nabla u_{\leq q-2}\,dx$. Therefore, we infer
\begin{equation}\notag
\begin{split}
|K_{12}|\leq &\sum_{q\geq -1}\lambda_q^{2s}\int_{\R^3}|(\nabla\Q_q\otimes\nabla\Q_q)\nabla u_{\leq q-2}|\,dx\\
\leq&\sum_{q\leq Q+2}\lambda_q^{2s}\int_{\R^3}|(\nabla\Q_q\otimes\nabla\Q_q)\nabla u_{\leq q-2}|\,dx\\
&+\sum_{q>Q+2}\lambda_q^{2s}\int_{\R^3}|(\nabla\Q_q\otimes\nabla\Q_q)\nabla u_{\leq Q}|\,dx\\
&+\sum_{q>Q+2}\lambda_q^{2s}\int_{\R^3}|(\nabla\Q_q\otimes\nabla\Q_q)\nabla u_{(Q,q-2]}|\,dx\\
\equiv&K_{121}+K_{122}+K_{123},
\end{split}
\end{equation}
with
\begin{equation}\notag
\begin{split}
|K_{121}|
\leq&\sum_{q\leq Q+2}\lambda_q^{2s}\int_{\R^3}|(\nabla\Q_q\otimes\nabla\Q_q)\nabla u_{\leq q-2}|\,dx\\
\lesssim &\sum_{q\leq Q+2}\lambda_q^{2s}\|\nabla u_{\leq q-2}\|_\infty\|\nabla \Q\|_2^2\\
\lesssim &Qf(t)\sum_{q\leq Q+2}\lambda_q^{2s+2}\|\Q\|_2^2;
\end{split}
\end{equation}
\begin{equation}\notag
\begin{split}
|K_{122}|
\leq&\sum_{q> Q+2}\lambda_q^{2s}\int_{\R^3}|(\nabla\Q_q\otimes\nabla\Q_q)\nabla u_{\leq Q}|\,dx\\
\lesssim &\sum_{q> Q+2}\lambda_q^{2s}\|\nabla u_{\leq Q}\|_\infty\|\nabla \Q\|_2^2\\
\lesssim &Qf(t)\sum_{q> Q+2}\lambda_q^{2s+2}\|\Q\|_2^2;
\end{split}
\end{equation}
and
\begin{equation}\notag
\begin{split}
|K_{123}|
\leq&\sum_{q> Q+2}\lambda_q^{2s}\int_{\R^3}|(\nabla\Q_q\otimes\nabla\Q_q)\nabla u_{(Q,q-2]}|\,dx\\
\lesssim &\sum_{q> Q+2}\lambda_q^{2s}\|\nabla u_{(Q,q-2]}\|_\infty\|\nabla \Q\|_2^2\\
\lesssim &\sum_{q> Q+2}\lambda_q^{2s+2}\|\Q\|_2^2\sum_{Q<p\leq q-2}\lambda_p\|u_p\|_\infty\\
\lesssim &c_r\mu\sum_{q> Q+2}\lambda_q^{2s+2}\|\Q\|_2^2\sum_{Q<p\leq q-2}\lambda_p^2\\
\lesssim &c_r\mu\sum_{q> Q+2}\lambda_q^{2s+4}\|\Q\|_2^2.
\end{split}
\end{equation}
To estimate $K_{3}$, we use integration by parts first, and then split it as  
\begin{equation}\notag
\begin{split}
|K_{3}|\leq&\left|\sum_{p\geq -1}\lambda_q^{2s}\sum_{-1\leq q\leq p+2}\int_{\mathbb R^3}\Delta_q(u_p\cdot\nabla) \tilde \Q_p\Delta\Q_q \, dx\right|\\
\leq&\sum_{p>Q}\lambda_q^{2s}\sum_{-1\leq q\leq p+2}\int_{\mathbb R^3}\left|\Delta_q(u_p\cdot\nabla)\tilde \Q_p\Delta\Q_q\right| \, dx\\
&+\sum_{-1\leq p\leq Q}\lambda_q^{2s}\sum_{-1\leq q\leq p+2}\int_{\mathbb R^3}\left|\Delta_q(u_p\cdot\nabla)\tilde \Q_p\Delta\Q_q\right| \, dx\\
\equiv & K_{31}+K_{32},
\end{split}
\end{equation}
with
\begin{equation}\notag
\begin{split}
|K_{31}|
\lesssim &\sum_{p> Q}\|u_p\|_\infty\|\nabla \Q_p\|_2\sum_{-1\leq q\leq p+2}\lambda_q^{2s}\|\Delta \Q_q\|_2\\
\lesssim &c_r\mu\sum_{p> Q}\lambda_p\|\nabla\Q_p\|_2\sum_{-1\leq q\leq p+2}\lambda_q^{2s+2}\|\Q_q\|_2\\
\lesssim &c_r\mu\sum_{p> Q}\lambda_p^{s+2}\|\Q_p\|_2\sum_{-1\leq q\leq p+2}\lambda_q^{s+2}\|\Q_q\|_2\lambda_{q-p}^{s}\\
\lesssim &c_r\mu\sum_{q\geq -1}\lambda_q^{2s+4}\|\Q_q\|_2^2;
\end{split}
\end{equation}
and 
\begin{equation}\notag
\begin{split}
|K_{32}|
\lesssim &\sum_{-1\leq p\leq Q}\|u_p\|_\infty\|\nabla\Q_p\|_2\sum_{-1\leq q\leq p+2}\lambda_q^{2s}\|\Delta\Q_q\|_2\\
\lesssim &f(t)\sum_{-1\leq p\leq Q}\|\Q_p\|_2\sum_{-1\leq q\leq p+2}\lambda_q^{2s+2}\|\Q_q\|_2\\
\lesssim &f(t)\sum_{-1\leq p\leq Q}\lambda_p^{s+1}\|\Q_p\|_2\sum_{-1\leq q\leq p+2}\lambda_q^{s+1}\|\Q_q\|_2\lambda_{q-p}^{s+1}\\
\lesssim &f(t)\sum_{-1\leq q\leq Q+2}\lambda_q^{2s+2}\|\Q_q\|_2^2.
\end{split}
\end{equation}
Following the above analysis and computation, we obtain
\begin{equation}\label{est-j2}
|J_2|+|K|\lesssim c_r\mu\sum_{q\geq -1}\lambda_q^{2s+4}\|\Q_q\|_2^2+Qf(t)\sum_{q\geq -1}\lambda_q^{2s+2}\|\Q_q\|_2^2.
\end{equation}

%(\ref{est-m}), (\ref{est-i}), (\ref{est-j1}),
Combining (\ref{ineq-energy})--(\ref{est-l}) and (\ref{est-j2}) yields that for some small enough constant $c_r$ with $2\leq r<\frac3s$ and $\frac12<s<1$
\begin{equation}\notag%\label{energy2}
\frac{d}{dt}\sum_{q\geq -1}\lambda_q^{2s}\left(\|u_q\|_2^2+\|\nabla\Q_q\|_2^2\right)
\lesssim \left(Q(t)f(t)+1\right)\sum_{q\geq -1}\lambda_q^{2s}\left(\|u_q\|_2^2+\|\nabla\Q_q\|_2^2\right),
\end{equation}
i.e., there exists an adimensional constant $C=C(r,\nu,\mu, s)$, such that 
\begin{equation}\label{energy3}
\frac{d}{dt}\left(\|u\|_{\dot H^s}^2+\|\nabla\Q\|_{\dot H^s}^2\right)
\leq C\left(Q(t)f(t)+1\right)\left(\|u\|_{\dot H^s}^2+\|\nabla\Q\|_{\dot H^s}^2\right) .
\end{equation}
We claim that, for $t>0$ 
\begin{equation}\label{Q2}
Q(t)\leq C(\nu,\mu,s)\left(1+\log {\|u(t)\|_{\dot H^s}}\right)
\end{equation}
Indeed, it follows from (\ref{Q}) and Bernstein's inequality that
\begin{equation}\notag
\Lambda(t)\leq (c_r\min\{\nu,\mu\})^{-1}\Lambda(t)^{\frac3r}\|u_{Q(t)}(t)\|_r 
\leq (c_r\min\{\nu,\mu\})^{-1}\Lambda(t)^{\frac32}\|u_{Q(t)}(t)\|_2. 
\end{equation}
Thus, we obtain
\[\Lambda^{s-\frac12}(t) \leq (c_r\min\{\nu,\mu\})^{-1}\|u(t)\|_{\dot H^s}. \]
Since $s>\frac12$, (\ref{Q2}) follows immediately. 
%Notice that this holds as long as  $r,s$ satisfy $1/2< s<1$ and $2\leq r<3/s$, and hence for any $r$ belonging to $[2, 6)$.
Combining (\ref{energy3}) and (\ref{Q2}) yields 
\begin{equation}\notag
\frac{d}{dt}\left(\|u\|_{\dot H^s}^2+\|\nabla\Q\|_{\dot H^s}^2\right)
\leq C(\nu,\mu, r,s)f(t)\left(1+\log {\|u(t)\|_{\dot H^s}}\right)\left(\|u\|_{\dot H^s}^2+\|\nabla\Q\|_{\dot H^s}^2\right).
\end{equation}
Therefore, due to the assumption $f\in L^1(0,T)$, applying Gr\"onwall's inequality to the above inequality gives us that $\|u(t)\|_{\dot H^s}^2+\|\nabla\Q(t)\|_{\dot H^s}^2$ is bounded on $[0,T)$. It concludes the proof of Theorem \ref{thm}.

%With more delicate analysis as in \cite{CDmhd}, one can prove Theorem \ref{thm-strong} which provides a weaker condition to guarantee global regularity. We omit the proof and leave it to interested readers. 

%\Endrefs
\end{document}